\documentclass[a4paper,12pt]{article}
\usepackage{latexsym}
\usepackage{amsfonts}
\title{An answer to a question of Pyrih\thanks{1991 Math.\
Subject Classification --- Primary: 54G05; Secondary: 54G99,
54H05. \protect\newline Key words and phrases --- extremally
disconnected, open-normal, discrete.}}
\author{Julian Dontchev\\Department of Mathematics\\University
of Helsinki\\PL 4, Yliopistonkatu 5\\00014 Helsinki\\Finland}
\date{}
\begin{document}
\baselineskip=20pt plus 1pt minus 1pt
\maketitle

A topological space $(X,\tau)$ is called {\em extremally
disconnected} \cite{St1} if the closure of every open set is
open. Recently, Pyrih \cite{P1} defined a space $X$ to be {\em
open-normal} if for any two disjoint open sets $A$ and $B$ there
exist disjoint closed sets $F_A$ and $F_B$ such that $A \subset
F_A$ and $B \subset F_B$. He asked if every open-normal space is
discrete. The answer to his question is no.

{\bf Proposition 1.} {\em A topological space $(X,\tau)$ is
open-normal if and only if it is extremally disconnected.}

{\em Proof.} Assume first that $X$ is open-normal. Let $A$ be an
open subset of $X$. Then, $A$ and $B = X \setminus {\rm cl} (A)$
are disjoint open sets and hence there exist disjoint closed
sets $F_A$ and $F_B$ such that $A \subset F_A$ and $B \subset
F_B$. Clearly, $F_A = {\rm cl} (A)$ and $F_B = B$. Hence, ${\rm
cl} (A) = X \setminus F_B$ is open. Thus, $X$ is extremally
disconnected.

Assume next that $X$ is extremally disconnected. If $A$ and $B$
are disjoint open subsets of $X$, then ${\rm cl} (A)$ and $X
\setminus {\rm cl} (A)$ are disjoint open sets containing $A$
and $B$ respectively. $\Box$

It is well-known that there exist extremally disconnected
non-discrete spaces: the Stone-\v{C}ech compactification of
every discrete space, etc. A detailed bibliography on the recent
progress of the study of extremally disconnected spaces may be
found in \cite{DR1}.

In \cite{P1}, Pyrih proved that every extremally disconnected
metric space is discrete. A much better result exists. In
\cite{G1}, Gleason proved that every convergent sequence
$(x_{1},x_{2},\ldots)$ of an extremally disconnected Hausdorff
space is stationary. Recall that a sequence $(x_1,x_2,\ldots)$
is called {\em stationary} if for some $n$, we have $x_n =
x_{n+1} = x_{n+2} = \ldots$. As a consequence of Gleason's
result we have the following:

{\bf Proposition 2.} {\em Every sequential Hausdorff space
which is extremally disconnected is discrete.}

Note that the following implications hold and none of them is
reversible:

Metric $\Rightarrow$ First countable and Hausdorff $\Rightarrow$
Fr\'{e}chet and Hausdorff $\Rightarrow$ Sequential and Hausdorff.

\
E-mail: {\tt dontchev@cc.helsinki.fi}
\

\baselineskip=12pt

\begin{thebibliography}{4}\frenchspacing


\bibitem{DR1} {J.~Dontchev and D.~Rose}, {Extremal
disconnectedness modulo dual filtrations}, {\em Math.
Pannonica}, {\bf 9} (1998), 33-46.

\bibitem{G1} {A.M.~Gleason}, {Projective topological spaces},
{\em Illinois J. Math.}, {\bf 2} (1958), 482--489.

\bibitem{P1} {P.~Pyrih}, {A space where disjoint open sets have
disjoint closures}, {\em Questions Answers Gen. Topology}, {\bf
16} (2) (1998), 133--134.

\bibitem{St1} {M.H.~Stone}, {Algebraic characterizations of
special Boolean rings}, {\em Fund. Math.}, {\bf 29} (1937),
223--302.


\end{thebibliography}
\end{document}